\documentclass[12pt]{article}
\usepackage{amssymb,amsmath,amsfonts,amssymb}
\usepackage{graphics,graphicx,color,hhline}
\usepackage{bbm} 
\usepackage{euscript}  
\usepackage{authblk}

 \usepackage{mathtools}

\def\@abssec#1{\vspace{.05in}\footnotesize \parindent .2in 
{\bf #1. }\ignorespaces} 

\graphicspath{{/EPSF/}{Figures/}}

\newtheorem{theorem}{Theorem}[section]
\newtheorem{lemma}[theorem]{Lemma}

\newtheorem{assumption}[theorem]{Assumptions}
\def \Rm {\mathbb R}
\def \NN {\mathbb N}

\newcommand{\rhoeq}[1]{\varrho_{\textrm{eq}}[#1]}
\newcommand{\eps}{\varepsilon}

\newcommand{\vt}{\vartheta}
\newcommand{\mn}{\mathrm{n}}

\newcommand{\ds}{\displaystyle}

\newcommand{\calQ}{\mathcal Q}

\newcommand{\calE}{\mathcal E}
\newcommand{\Tr}{\textnormal{Tr}}

\newcommand{\calJ}{\mathcal J}

\def\fref#1{{\rm (\ref{#1})}}

\newcommand{\cout}[1]{}

\newcommand{\be}{\begin{equation}}
\newcommand{\ee}{\end{equation}}
\newcommand{\bea}{\begin{eqnarray}}
\newcommand{\eea}{\end{eqnarray}}
\newcommand{\bee}{\begin{eqnarray*}}
\newcommand{\eee}{\end{eqnarray*}}
\newcommand{\bal}{\begin{align*}}
                    \newcommand{\eal}{\end{align*}}

\def\H2p{{H}^2}

\DeclarePairedDelimiter\bra{\langle}{\rvert}
\DeclarePairedDelimiter\ket{\lvert}{\rangle}

\begin{document}
{\title{An asymptotic preserving scheme for the quantum Liouville-BGK equation}}

 \author{Romain  Duboscq \footnote{Romain.Duboscq@math.univ-tlse.fr}}
 \affil{Institut de Math\'ematiques de Toulouse ; UMR5219\\Universit\'e de Toulouse ; CNRS\\INSA, F-31077 Toulouse, France}
 \author{Olivier Pinaud \footnote{olivier.pinaud@colostate.edu}}
 \affil{Department of Mathematics, Colorado State University\\ Fort Collins CO, 80523}

\maketitle

\begin{abstract}
We are interested in this work in the numerical resolution of the Quantum Liouville-BGK equation, which arises in the derivation of quantum hydrodynamical models from first principles. Such models are often obtained in some asymptotic limits, for instance a diffusion or a fluid limit, and as a consequence the original Liouville equation contains small parameters. A standard method such as a split-step algorithm is then accurate provided the time step is sufficiently small compared to the asymptotic parameter, which is a severe limitation. In the case of the diffusion limit, we propose a numerical method that is accurate for time steps independent of the small parameter, and which captures well both the microscopic dynamics and the diffusion limit. Our approach is substantiated by an informal theoretical error analysis.
\end{abstract}

\section{Introduction}
This work is concerned with the discretization of the quantum Liouville-BGK equation (QLE) of the form, after rescaling,
\begin{equation} \label{QLEintro}
    i\eps \partial_t\varrho = [H,\varrho] + \frac{i}{\eps}Q(\varrho),
\end{equation}
where $\eps \ll 1$, $\varrho$ is a density operator, i.e. a trace class self-adjoint nonnegative operator on some Hilbert space, $[\cdot,\cdot]$ denotes the commutator between two operators, $H$ is a given Hamiltonian, and $Q$ is a BGK-type collision operator  \cite{BGK} of the form
\begin{equation*}
    Q(\varrho) = \rhoeq{\varrho}-\varrho.
\end{equation*}
Above, $\rhoeq{\varrho}$ is a quantum statistical equilibrium obtained by minimizing the quantum free energy $F$, defined by, for appropriate density operators $\sigma$, 
\begin{equation}
\label{FE}
    F(\sigma) = T_e\, \Tr(\sigma \log \sigma - \sigma) + \Tr(H \sigma),
\end{equation}
under the constraint that the local density of particles of $\sigma$ is the same as that of $\varrho$, with $\varrho\equiv \varrho(t)$ the solution to \fref{QLEintro}. In other words, if $W_\sigma$ and $W_\varrho$ are the Wigner transforms of $\sigma$ and $\varrho$, then this constraint is expressed mathematically as
$$n_\sigma(x):=\int_{\Rm^d} W_\sigma(x,p) dp=\int_{\Rm^d}  W_\varrho(x,p) dp=n_\varrho(x),$$
(further we will use more a convenient form for the definition of $n_\sigma$, see \fref{eq:moments}). In \fref{FE}, $\Tr(\cdot)$ denotes operator trace, $T_e$ is the temperature, and if $\sigma$ has integral kernel $\sigma(x,y)$, its Wigner transform is defined by \cite{LP}
$$
W_\sigma(x,p)=\frac{1}{(2 \pi)^d}\int_{\Rm^d} e^{i y \cdot p} \sigma(x-y/2,x+y/2) dy.
$$

The nondimensional parameter $\eps$ is the ratio of the so-called mean free time, which quantifies the average duration between collision events, and a typical simulation time. The limit $\eps \to 0$ is commonly referred to as the diffusion limit. This problem is motivated by a series of papers by Degond and Ringhofer on the derivation of quantum hydrodynamical models from first principles. We refer the reader to \cite{DR,DGMRlivre,QHD-CMS,isotherme,QHD-review,QET,jungel-matthes-milisic,jungel-matthes,jungelbook} for additional motivations and references. A critical difficulty in the resolution of \fref{QLEintro} is the presence of the nonlinear term $Q(\varrho)$, and a numerical scheme was proposed in \cite{BragP} in the case $\eps=1$. The method is based on a splitting technique, which reduces the treatment of the nonlinear part to the resolution of a linear problem. A serious drawback of this method is that the time stepsize has to decrease to zero as $\eps \to 0$ for reasonable accuracy and consequently the computational cost becomes prohibitive. Our main motivation in the present work is to propose an alternative method in which the time stepsize can be chosen independently of $\eps$, which then reduces the computational burden. Such techniques belong to the realm of asymptotic preserving methods (AP), see e.g. \cite{JinAP,lemouAP,jin2022asymptotic,liu2010analysis}.

In the kinetic case where \fref{QLEintro} is a classical transport equation, asymptotic preserving schemes are typically obtained by perturbation around the equilibrium and by projecting the kinetic equation onto the orthogonal complement of the kernel of $Q$, see e.g. \cite{lemouAP}. While such a procedure is quite straightforward in the kinetic case as the collision operator only acts on the velocities defined in phase space, it is unclear how this approach can be transposed to the quantum case where $Q$ is non-local and velocities do not form their own set of independent variables.

We then follow a different route based on the asymptotic expansion of $\varrho$ obtained when deriving the diffusion limit. As $\eps \to 0$, $\varrho$ converges to a so-called quantum Maxwellian, and its local density $n_\varrho$ solves the Quantum Drift-Diffusion (QDD) equation which plays a crucial role for the design of the AP scheme. More precisely, it is shown in \cite{QDD-JCP} that,  as $\eps \to 0$, $\varrho(t)$ converges to
$$\varrho_{\mathrm{eq}}[\varrho(t)] =e^{-\frac{1}{T_e}(H + \mathsf{A}[\mn(t)])},$$ where $\mathsf{A}[\mn(t)](x)$ is the (real-valued) chemical potential associated with the density constraint $\mn(t)=n_{\varrho(t)}$. The expression above is meant in the functional calculus sense. The potential $\mathsf{A}[\mn(t)](x)$ is such that the local density of $\varrho_{\mathrm{eq}}[\varrho(t)]$ is equal to $\mn(t)$, where $\mn(t)$ verifies the QDD equation
\begin{equation}\label{eq:qddint}
 \partial_t \mathrm{n}(t)  + \nabla \cdot \left(\mathrm{n}(t) \nabla \mathsf{A}[\mathrm{n}(t)]\right)=0.
\end{equation}
 
Our method to obtain an AP scheme consists of two steps. In the first one, we compute the local density of $\varrho(t)$ solution to the QLE equation using a modified version of \fref{eq:qddint} with a non-zero right hand side that depends on $\eps$. We adapt for this the techniques of \cite{QDD-SIAM} based on the minimization of a functional, and the important point is that the resolution can be done with a time stepsize independent of $\eps$. Once the local density is known, the second step consists in calculating the density operator $\varrho(t)$ using an asymptotic formula that is uniform in the time stepsize and in $\eps$. The resulting scheme has the AP property in that it is accurate independently of the value of $\eps$, and we have a diagram of the form
\begin{eqnarray} \label{diag}
\begin{array}{cccc}
\varrho^{\varepsilon}_{\delta t} & \stackrel{\varepsilon\to 0}{\longrightarrow} & \varrho_{\delta t}& \\
 \downarrow & \hspace{0.5em}\!\!\!\!\!\!\!\!\!\!\!\!\!\!\!\!\!\stackrel{\delta t\to 0}{}  & \downarrow
 &\!\!\!\!\!\!\! \stackrel{\delta t\to 0}{}  \\
\varrho^{\varepsilon} &\stackrel{\varepsilon\to 0}{\longrightarrow }  &\varrho &
\end{array}
\end{eqnarray}
Above, $\varrho^\eps$ denotes the solution to \fref{QLEintro} for fixed $\eps \in (0,1]$, $\varrho^\eps_{\delta t}$ is the numerical solution of our AP scheme for a time step $\delta t$, $\varrho$ is the exact diffusion limit, and $\varrho_{\delta t}$ is the diffusion limit of the AP scheme. Since the key issue is the treatment of the time variable and not of the spatial variables, our strategy is applicable to spatial domains of any dimensions and our method will be presented for an arbitrary spatial domain. 

Our main goal in this work is the derivation the AP scheme. A mathematically rigorous analysis of the  scheme is not possible at the present time due to the lack of a rigorous mathematical theory for QLE and QDD. The question is indeed quite difficult and hinges on the properties of the nonlinear map $\varrho \mapsto \varrho_{\mathrm{eq}}[\varrho]$ for which little is actually known. To the best of our knowledge, one can find existence results for \fref{QLEintro} and \fref{eq:qddint} in \cite{pinaudPoincare,MP-JDE, matthes}, while the minimization of the free energy under constraints is adressed in \cite{DP-MB,DP-JMPA,DP-JMP,DP-CVPDE,DP-JFA,MP-KRM,MP-JSP}, but a general theory is elusive. This does not prevent us from deriving the AP scheme, and our analysis will then remain mostly informal. We will hence study the scheme under some realistic assumptions on the solutions to \fref{eq:qddint} and on the map $\varrho \mapsto \varrho_{\mathrm{eq}}[\varrho]$. Under such conditions, we establish uniform error estimates backing our claim that the scheme is AP. Numerical simulations will be the object of a future work.

The paper is structured as follows: Section \ref{models} is dedicated to the derivation of the AP scheme and to some of its properties; in Section \ref{proofs}, we perform an informal theoretical error backing the claim that the numerical scheme is asymptotic preserving.

\paragraph{Acknowledgment.} This work was supported by NSF grant DMS-2006416 and DMS-2404785.

\section{The AP scheme} \label{models}

In this Section, we introduce some notation, present our methodology, explicitate the AP scheme and discuss some of its properties.

\paragraph{Setting of the problem.} 

We first write a density operator $\varrho$ in terms of its spectral elements,
\[\varrho = \sum_{p\in \mathbb{N}} \rho_{p}\ket{\psi_{p}}\bra{\psi_{p}},\]
where we used the Dirac bra-ket notation, and where $\{\rho_p,\psi_p\}$ are the $p-$th eigenvalue and eigenfunction pair for $\varrho$, eigenvalues counted with multiplicity. In our problem of interest, the density operators are typically full-rank, that is all eigenvalues are strictly positive, and form then a sequence $\{\rho_p\}_{p \in \NN}$ decreasing to zero. This is a consequence of the fact, proved in \cite{MP-JSP}, that the equilibrium $\rhoeq{\varrho}$ is full-rank. With this notation, the local density $n_\varrho$ and local current $j_\varrho$ associated to $\varrho$ are defined by (we will use both notations $n_\varrho$ and $n[\varrho]$, $j_\varrho$ and $j[\varrho]$, depending on which one is more convenient)
\begin{equation}\label{eq:moments}
n_\varrho(x) = \sum_{p\in\NN} \lambda_p |\psi_p(x)|^2\quad\mathrm{and}\quad j_\varrho(x) = 2\Im\left(\sum_{p\in\NN}\lambda_p \psi_p^*(x)\nabla \psi_p(x)  \right).
\end{equation}
The local density can also be equivalently defined by duality in terms of the trace operator  $\Tr(\cdot)$, i.e., if $\Omega \subseteq \Rm^{d}$, $d \geq 1$, is the spatial domain, 
 \[(n_{\varrho},\psi) := \int_{\Omega} n_{\varrho}\psi dx = \Tr(\varrho\psi),\] for all smooth function $\psi$ (we identify $\psi$ with the corresponding multiplication operator).

In the context of particle transport in nanostructures, the Hamiltonian $H$ in \fref{QLEintro} is of the form 
\[  H = H_0+V,  \quad \text{where } \quad H_0= -\Delta,\]
equipped with appropriate boundary conditions that ensure that $H$ and $H_0$ are self-adjoint, and where $V$ is a bounded real-valued potential. While we only consider a linear Hamiltonian in this work, our method can straightforwardly be generalized to the case of nonlinear potentials of Poisson type treated in \cite{BragP}.

With $L(\varrho) : = [H,\varrho]$, we rewrite the QLE as
\begin{equation}\label{eq:Liouville-BGK}
i\varepsilon\partial_t \varrho^{\varepsilon}(t) = L(\varrho^{\varepsilon}(t)) + \frac{i}{\varepsilon} Q(\varrho^{\varepsilon}(t)),
\end{equation}
with $\varrho^{\varepsilon}(0) = \varrho_0$. We will write $\mathrm{n}^{\varepsilon}(t) := n[\varrho^{\varepsilon}(t)]$ for simplicity. Since $\varrho_{\mathrm{eq}}[\varrho^\eps(t)]$ only depends on $\varrho^\eps(t)$ through $\mathrm{n}^{\varepsilon}(t)$, we denote it by $\vt[\mathrm{n}^{\varepsilon}](t)$.

We recall that the equilibrium in $Q$ is obtained by minimizing \fref{FE} under the constraint $n_\sigma=n_{\varrho^\eps(t)}$. For a density $n_0(x)$, it is shown e.g. in \cite{QDD-JCP} that the equilibrium has the form, in the functional calculus sense, 
\begin{equation}\label{eq:exprminimiser}
\vartheta[n_0] =  e^{-\frac{1}{T_e}(H_0 + \mathsf{A}[n_0])},
\end{equation}
where $\mathsf{A}[n_0](x)$ the chemical potential associated with the density constraint. We recall that the properties of the map $n_0 \mapsto \vartheta[n_0]$, in particular its Lipschitz character, are largely an open problem. Yet, enough is known to show that \fref{QLEintro} admits a unique solution in 1D in the proper functional setting \cite{MP-JDE,BragP}.

\paragraph{General strategy.}
The starting point is to introduce the operator
\begin{equation*}
\sigma^{\varepsilon}(t) := \frac{1}{\varepsilon}Q(\varrho^{\varepsilon}(t)) - i L(\vartheta[\mathrm{n}^{\varepsilon}](t)) = \frac{1}{\varepsilon}(\vartheta[\mathrm{n}^{\varepsilon}](t) - \varrho^{\varepsilon}(t)) - i L(\vartheta[\mn^{\varepsilon}](t)),
\end{equation*}
so that
\begin{equation} \label{eqrho}
\varrho^{\varepsilon}(t) = \vartheta[\mathrm{n}^{\varepsilon}](t) -  \varepsilon \sigma^{\varepsilon}(t) - i \varepsilon L(\vartheta[\mathrm{n}^{\varepsilon}](t)).
\end{equation}
By construction, the collision operator $Q$ preserves the local density, which is equivalent to $n[Q(\varrho^{\varepsilon}(t))]=0$. The first-order continuity relation for QLE then yields
  \be \label{conti}
\varepsilon \partial_t \mn^\varepsilon(t) +\nabla \cdot j[\varrho^{\varepsilon}(t)]=0.
\ee
We now explicitate $j[\varrho^{\varepsilon}(t)]$. We remark first that the equilibrium $\vt[\mn^\eps]$ carries no current, and therefore $j[\vt[\mn^\eps]]=0$. This is easily seen, for instance, by choosing real-valued eigenfunctions for $H_0+\mathsf{A}[n^\eps]$, where $\mathsf{A}$ is the chemical potential in expression \eqref{eq:exprminimiser}. Moreover, it is shown in \cite{QDD-JCP} in the derivation of QDD that $j[iL(\vartheta[\mathrm{n}^{\varepsilon}])] = -\mathrm{n}^{\varepsilon} \nabla \mathsf{A}[\mathrm{n}^{\varepsilon}]$. This then yields that
$$
j[\varrho^{\varepsilon}(t)]=\eps \mathrm{n}^{\varepsilon} \nabla \mathsf{A}[\mathrm{n}^{\varepsilon}(t)]-\eps j[\sigma^{\varepsilon}(t)].
$$
Injecting the latter into \fref{conti}, we obtain the equation
\begin{equation}\label{eq:neq}
 \partial_t \mathrm{n}^{\varepsilon}(t)  + \nabla \cdot \left(\mathrm{n}^{\varepsilon}(t) \nabla \mathsf{A}[\mathrm{n}^{\varepsilon}(t)]\right)=  \nabla \cdot j[\sigma^{\varepsilon}(t)],
\end{equation}
which is exact at this point. Note that it is shown in \cite{QDD-JCP} that $j[\sigma^{\varepsilon}] \to 0$ as $\eps \to 0$, leading to the QDD equation \fref{eq:qddint}.  The key to the AP scheme is then to derive an approximation of $\sigma^{\varepsilon}$ that is uniform in $\eps$ and yields the exact value as the time stepsize goes to zero. Once this approximation is available, one can solve \fref{eq:neq} from $t$ to $t+\delta t$ using similar methods as those for QDD with a time stepsize $\delta t$ that is independent of $\eps$. One can then iterate: when the density $\mn^\eps(t+\delta t)$ is known, it is injected in the approximate expression of $\sigma^\eps$ that can be evaluated numerically uniformly in $\eps$. The density operator $\varrho^\eps(t)$ follows from expression \fref{eqrho}.

\paragraph{The AP scheme.} We now derive the scheme explicitly. For $\tau>0$ and a self-adjoint operator $\sigma$, let
\begin{equation} \label{defkap}
\kappa_{\varepsilon}(\tau) : = 1-e^{-\varepsilon^{-2}\tau} - \varepsilon^{-2}\tau e^{-\varepsilon^{-2}\tau}\quad\mbox{and}\quad
S_{\varepsilon,\tau}(\sigma) : = e^{-\tau}e^{-i\varepsilon\tau H }\sigma e^{i\varepsilon\tau H}.
\end{equation}
The operator $S_{\varepsilon,\tau}(\sigma)$ is simply the (exponentially damped) free evolution of QLE with $\rhoeq{\varrho^\eps}=0$ and initial condition $\sigma$. We have then the following crucial lemma:
\begin{lemma}\label{lem:expsigma}
For any $t> s\geq 0$, the following equality holds
\begin{align}
  -\varepsilon \sigma^{\varepsilon}(t) =&  S_{\varepsilon,(t-s)/\varepsilon^2}[\varrho^\eps(s)] -\vartheta[\mathrm{n}^{\varepsilon}](t)e^{-(t-s)/\eps^2} \nonumber
  + i\varepsilon [H,\vartheta[\mathrm{n}^{\varepsilon}](t) ]\left(1-\kappa_{\varepsilon}(t-s)\right)
 \\&+ \varsigma^{\varepsilon}_1(s,t),\label{eq:expsigma1}
\end{align}
where $\varsigma^{\varepsilon}_1=O(\eps^2 \wedge (t-s) )$ and $\nabla \cdot j[\varsigma^{\varepsilon}_1]=O(\eps [\eps^2 \wedge (t-s)])$.
\end{lemma}

The proof of the Lemma is given in Section \ref{proofs}. Note that the error term $\varsigma^{\varepsilon}_1$ is uniform both in $\eps$ and in $(t-s)$. Neglecting it in expression \fref{eq:expsigma1} then yields a uniform error in $\eps$ and in $(t-s)$, which we exploit in the construction in the AP scheme.

Using once more the facts that $j[ \vartheta[\mathrm{n}^{\varepsilon}]]=0$ and that $$j[iL(\vartheta[\mathrm{n}^{\varepsilon}])]= j[i [H,\vartheta[\mathrm{n}^{\varepsilon}])]= -\mathrm{n}^{\varepsilon} \nabla \mathsf{A}[\mathrm{n}^{\varepsilon}],$$
we find from the above Lemma that
$$
j[\sigma^{\varepsilon}(t)] \simeq  -\eps^{-1}j[S_{\varepsilon,(t-s)/\varepsilon^2}[\varrho^\eps(s)]]+\left(1-\kappa_{\varepsilon}(t-s)\right)\mathrm{n}^{\varepsilon} \nabla \mathsf{A}[\mathrm{n}^{\varepsilon}].
$$
This gives the approximate relation, for $t \geq s$,
\begin{equation}\label{eq:neq2}
 \partial_t \mathrm{n}^{\varepsilon}(t)  + \kappa_{\varepsilon}(t-s)\nabla \cdot \left(\mathrm{n}^{\varepsilon}(t) \nabla \mathsf{A}[\mathrm{n}^{\varepsilon}(t)]\right)\simeq   - \eps^{-1} \nabla \cdot j[S_{\varepsilon,(t-s)/\varepsilon^2}[\varrho^\eps(s)]],
\end{equation}
where the error is uniform in $\eps$ and $(t-s)$. We will refer to the equation above as \textit{the modified QDD equation}, and it is probably best seen as a nonlinear evolution equation on $\mathsf{A}$ rather than on the density $\mn^\eps$. The term $\varrho^\eps(s)$ in \fref{eq:neq2} is known from the previous iteration, and as a consequence the right hand side is available provided we compute the free evolution $S_{\varepsilon,(t-s)/\varepsilon^2}[\varrho^\eps(s)]$. We then solve  \fref{eq:neq2} to obtain the density $\mn^\eps$ at time $t+\delta t$. Finally, we plug $\mn^\eps(t+\delta t)$ in \fref{eqrho}-\fref{eq:expsigma1} to obtain $\varrho^\eps(t+\delta t)$.

With $t_n = n\delta t$ and $\varrho^\eps_n=\varrho^\eps(t_n)$, an intermediate scheme where $\varsigma_1^\eps$ is neglected reads, for $n\geq 0$,
\begin{equation}
\left\{\begin{array}{ll}
\ds & \ds \tilde \mn^{\varepsilon}(t_{n+1}) -  \tilde \mn^{\varepsilon}(t_n)+  \int_{t_n}^{t_{n+1}} \kappa_{\varepsilon}(r-t_n)\nabla \cdot \left(\tilde \mn^{\varepsilon}(r) \nabla \mathsf{A}[\tilde \mn^{\varepsilon}(r)]\right) dr 
\ds \\ &\hspace{8em}\ds =  -\varepsilon^{-1} \int_{t_n}^{t_{n+1}} \nabla \cdot j[S_{\varepsilon,(r-t_n)/\varepsilon^2}[\tilde \varrho^{\varepsilon}_n]] dr
\ds \\ [3mm] &\tilde \varrho^{\varepsilon}_{n+1} = S_{\varepsilon,\delta t/\varepsilon^2}[\tilde \varrho^{\varepsilon}_n] +  \vartheta[\tilde \mn^{\varepsilon}](t_{n+1})(1-e^{-\varepsilon^{-2}\delta t}) - \varepsilon \kappa_{\varepsilon}(\delta t) [iH,\vartheta[\tilde \mn^{\varepsilon}](t_{n+1}) ].
\end{array}\right.\label{eq:scheme1}
\end{equation}

It remains to discretize the propagator $e^{-i \tau H}$ and the integrals w.r.t. $r$ above. The function $\kappa_\eps$ in the l.h.s. and the real exponential term in $S_{\varepsilon,(r-t_n)/\varepsilon^2}$ in the r.h.s. both evolve at the scale $\eps^{-2}$ and will inevitably introduce large errors if discretized. We then integrate them exactly and proceed as follows. Following \cite{QDD-SIAM}, we use a first-order implicit discretization of the nonlinear term in the first equation in \fref{eq:scheme1}, which is then approximated by 
$$
\int_{t_n}^{t_{n+1}} \kappa_{\varepsilon}(r-t_n)\nabla \cdot \left(\tilde \mn^{\varepsilon}(r) \nabla \mathsf{A}[\tilde \mn^{\varepsilon}(r)]\right) dr \simeq \nabla \cdot \left(\tilde \mn^{\varepsilon}(t_n) \nabla \mathsf{A}[\tilde \mn^{\varepsilon}(t_{n+1})]\right) \int_{t_n}^{t_{n+1}} \kappa_{\varepsilon}(r-t_n) dr.
$$
Regarding the term $S_{\varepsilon,(r-t_n)/\varepsilon^2}[\tilde \varrho^{\varepsilon}_n]$, we write
$$
\int_{t_n}^{t_{n+1}} \nabla \cdot j[S_{\varepsilon,(r-t_n)/\varepsilon^2}[\tilde \varrho^{\varepsilon}_n]] dr=\int_{t_n}^{t_{n+1}} e^{-\eps^{-2}(r-t_n)}  f(r-t_n)dr
$$
where $f(r-t_n)=\nabla \cdot j[e^{-i (r-t_n) H/\eps}\tilde \varrho^{\varepsilon} e^{i (r-t_n) H/\eps}]$, and use a midpoint rule adapted to the exponential weight to discretize the integral. With
$$
a_\eps=\frac{\int_0^{\delta t} r e^{-r/\eps^2} dr}{\int_0^{\delta t} r e^{-r/\eps^2} dr} \in (0,\delta t),
$$
we find
\bee
\int_{t_n}^{t_{n+1}} \nabla \cdot j[S_{\varepsilon,(r-t_n)/\varepsilon^2}[\tilde \varrho^{\varepsilon}_n]] dr &\simeq&f(a_\eps) \int_{t_n}^{t_{n+1}} e^{-\eps^{-2}(r-t_n)}  dr\\
&\simeq&\eps^2 f(a_\eps) (1-e^{-\delta t/\eps^2}),
\eee
at third order in $\delta t$.

The final step is to approximate the free propagator in $S_{\varepsilon,\delta t/\varepsilon^2}[\tilde \varrho^{\varepsilon}_n]$ for the calculation of $\tilde \varrho_{n+1}^\eps$. For this, we use a globally second order time discretization, for instance a Crank-Nicolson scheme. We denote the approximation by $\hat S_{\varepsilon,\delta t/\varepsilon^2} \sigma =e^{-\delta t/\eps^2} U^\eps_{\delta_t} \sigma (U^\eps_{\delta_t})^*$ where $U^\eps_{\delta_t}$ is the discretization of $e^{-i \delta t H /\eps}$. Note that we are not able to establish theoretically the AP character of the scheme when the approximation of the propagator is only first order in time. 
With
$$
\xi_\eps(\delta t)=\int_0^{\delta t} \kappa_\eps(r)dr,
$$
the full semi-discrete scheme is
\begin{equation}
\left\{\begin{array}{ll}
         \ds & \ds \mn^{\varepsilon}_{n+1} -  \mn^{\varepsilon}_n+\xi_{\varepsilon}(\delta t)\nabla \cdot \left(\mn^{\varepsilon}_n \nabla \mathsf{A}[\mn_{n+1}^\eps]\right)= -\eps (1-e^{-\frac{\delta t}{\eps^2}}) \nabla \cdot j[U^\eps_{a_\eps} \varrho_n^\eps (U^\eps_{a_\eps})^*]
\ds \\ [3mm] & \varrho^{\varepsilon}_{n+1} = \hat S_{\varepsilon,\delta t/\varepsilon^2}[ \varrho^{\varepsilon}_n] +  \vartheta[ \mn^{\varepsilon}_{n+1}](1-e^{-\frac{\delta t}{\eps^2}}) - \varepsilon \kappa_{\varepsilon}(\delta t) [i H,\vartheta[ \mn^{\varepsilon}_{n+1}]].
\end{array}\right.\label{dQDD}
\end{equation}
The above system is the main result of this paper.
The first equation above can be solved by adapting the minimization procedure of \cite{QDD-SIAM}, this is detailed further. 
%

Under reasonable assumptions on the well-posedness and stability of the QLE and QDD equations, and on the stability of the map $n_0 \mapsto \vartheta[n_0]$,  we show informally in Section \ref{proofs} that the error in the approximation of the density $\mn^\eps(t_n)$ by $\mn^\eps_n$ is globally of order $\delta t$, uniformly in $\eps$ and $\delta t$. In terms of the density operator $\varrho^\eps$, this translates into a global error of order $\delta t$ as well. This shows the AP character of the scheme. While our error analysis is largely formal, it is still a helpful guide in the design of the AP scheme, in particular in the discretization of the time integrals in \fref{eq:scheme1}. 

\paragraph{Numerical methods.} We describe here how the semi-discrete modified QDD equation \fref{dQDD} can be solved numerically. We remark first that it is established in \cite{QDD-SIAM}, when $\kappa_\eps=1$ and the right-hand side vanishes, that \fref{dQDD} admits a unique solution for all $n \geq 0$ under appropriate conditions on the initial condition. The analysis in \cite{QDD-SIAM} directly generalizes to the modified equation. The solution is obtained by a minimization procedure as follows: denote by $\Omega$ the spatial domain, and equip \fref{dQDD} with either Neumann (as in \cite{QDD-SIAM}) or periodic boundary conditions. Adapting the method of \cite{QDD-SIAM}, the chemical potential $\mathsf{A}[\mn^\eps_{n+1}]$ at step $n+1$ is obtained by minimizing the functional
\be \label{J}
J(A)=\frac{\xi_\eps(\delta t)}{2}\int_{\Omega} \mn^\eps_n |\nabla A|^2 dx+\int_\Omega (\mn^\eps_n+f_n) Adx+ T_e \Tr \left(e^{-\frac{H+A}{T_e}}\right),
\ee
where $f_n(x)=-\eps (1-e^{-\frac{\delta t}{\eps^2}}) \nabla \cdot j[U^\eps_{a_\eps} \varrho_n^\eps (U^\eps_{a_\eps})^*]$, and where $\mn^\eps_n$, $\varrho_n^\eps$ are known. This functional is strictly convex and therefore admits a unique minimizer, whose associated Euler-Lagrange equation is the first equation in \fref{dQDD}. The functional can be minimized for instance by using a nonlinear gradient descent. 


For an initial condition $\varphi$, the calculation of the evolution $e^{-i \delta t  H/\eps} \varphi$ can be done with standard globally second order methods  for the resolution of the Schr\"odinger equation, such as a Crank-Nicolson or a Strang splitting scheme. 

\paragraph{Computational cost.} For one iteration, our AP scheme requires (i) the calculation of $e^{-i \delta t  H/\eps}\varphi$, which is needed in the free evolution of the density operator in \fref{dQDD}, and (ii) the minimization of the functional $J(A)$ defined in \fref{J}. This is similar in cost to what is calculated for the splitting scheme of \cite{BragP}, with the difference that in \cite{BragP} it is the free energy functional \fref{FE} that is minimized. There are additional commutators that have to be computed for the AP scheme but the cost is negligible compared to items (i) and (ii). So overall, the cost of the AP scheme for one iteration is comparable to that of the non-AP splitting scheme of \cite{BragP}. What makes the AP scheme nore advantageous is the fact that less iterations are needed since the time step $\delta t$ can be chosen independently of $\eps$.
\paragraph{Some properties the AP scheme.} We note first that the density $\mn^\eps$ obtained from \fref{dQDD} is positive. Indeed, the minimization problem associated to \fref{dQDD} returns $\mathsf{A}[\mn^{\varepsilon}_{n+1}]$, which is used to compute the equilibrium $e^{-\frac{1}{T_e}(H+\mathsf{A}[\mn^{\varepsilon}_{n+1}])}$ whose local density $\mn^\varepsilon_{n+1}$ is positive by construction. We remark in addition that the $\varrho_{n+1}^\eps$ obtained in the second equation of \fref{dQDD} is self-adjoint provided the initial density operator is self-adjoint. 

Regarding the positivity of the  operator $\varrho_{n+1}^\eps$ solution to the second equation in \fref{dQDD}, the question is to assess whether or not the operator
  \be \label{OP} \vartheta[\mn^{\varepsilon}](t)\left(1-e^{-\varepsilon^{-2}(t-s)}\right) - i\varepsilon [H,\vartheta[\mn^{\varepsilon}](t) ] \kappa_\eps(t-s)\ee
 is positive. Indeed, the first operator in the r.h.s. in the  second equation of \fref{eq:scheme1} is positive if the previous iterate is positive and we thus need to discuss the second term. As will be seen in Section \ref{proofs}, this operator is an approximation of the second term in the mild formulation \fref{mild}, and we have
  \begin{align*}
    \vartheta[\mn^{\varepsilon}](t)\left(1-e^{-\varepsilon^{-2}(t-s)}\right) - &i\varepsilon [H,\vartheta[\mn^{\varepsilon}](t) ] \kappa_\eps(t-s)\\
    &= \int_0^{\varepsilon^{-2}(t-s)}S_{\varepsilon,\tau}[\vartheta[\mn^{\varepsilon}](t-\varepsilon^{2}\tau)] d\tau-\varsigma_1^\eps(t,s).
  \end{align*}
  Let us denote by $T_\eps$ the first term on the r.h.s. above. It is positive since $\vartheta[\mn^{\varepsilon}](r)$ is positive, and we have already discussed the fact that $\vartheta[\mn^{\varepsilon}]$ is full-rank, namely that the smallest eigenvalue of $\vartheta[\mn^{\varepsilon}]$ is strictly positive. Call it $\lambda_0(r)$ and suppose $\lambda_m=\inf_{r \in [s,t]} \lambda_0(r)>0$, which is expected since $\vartheta[\mn^{\varepsilon}](r)$ remains full-rank at all times. Then, a short calculation shows that, for any $\psi$ with norm one, 
  $$
  (\psi,T_\eps \psi) \geq \lambda_m (1-e^{-(t-s)/\eps^2}).
  $$
  We have seen in Lemma \ref{lem:expsigma} that $\varsigma_1^\eps(t,s)$ is of order $\eps^2 \wedge (t-s)$.  When $\eps^2 \wedge (t-s)$ is sufficiently small, this implies then that the operator in \fref{OP} is positive and therefore so is $\varrho_{n+1}^\eps$. Note that there is also the possibility that the operator in \fref{OP} is negative but compensated by the first operator in the r.h.s. in the second equation in \fref{dQDD}, but this scenario is more difficult to analyze.

  \paragraph{Diffusion limit.} Sending $\eps$ to zero in \fref{dQDD} leads informally to the system

  $$
\left\{\begin{array}{ll}
         \ds & \ds \mn_{n+1} -  \mn_n+\delta t\nabla \cdot \left(\mn_n \nabla \mathsf{A}[\mn_{n+1}]\right)= 0
\ds \\ [3mm] & \varrho_{n+1} = \vartheta[ \mn^{\varepsilon}_{n+1}],
\end{array}\right.
$$
which is the implicit discretization of the QDD equation \fref{eq:qddint} of \cite{QDD-SIAM}. This is a further justification of the AP character of the scheme and of the diagram \fref{diag}.


\section{Theoretical analysis} \label{proofs}

We justify in this Section the AP character of the scheme. As explained in the Introduction, our analysis is mostly informal due to the lack of rigorous theories for QLE and QDD. Our main goal will be to derive error estimates on the density and density operator that are uniform in $\eps$ under reasonable assumptions on the dynamics. \medskip

We start with some preliminaries.


\paragraph{Preliminaries.} For a trace class  operator $\sigma$ on $L^2(\Omega)$ (we recall that $\Omega$ is the spatial domain), we denote by $\| \sigma \|_{\calJ_1}$ the trace norm, that is $\| \sigma \|_{\calJ_1}=\Tr |\sigma|$ where $|\sigma|=\sqrt{\sigma^* \sigma}$ for $\sigma^*$ the adjoint of $\sigma$. With $H_0=-\Delta$, equipped with boundary conditions on $\partial \Omega$ that turn $H_0$ into a self-adjoint operator (e.g. Neumann or periodic boundary conditions when $\Omega$ is bounded), we set
$$
\|\sigma\|_{\calE^2}=\Tr |\sigma| + \Tr (H_0 |\sigma| H_0).
$$
With $H=H_0+V$ for some bounded real-valued potential $V$, we will use the relation
\be \label{commu}
j[[H,[H,\vartheta[n_0]]]]=0,
\ee
where $\vartheta[n]$ is the equilibrium associated with a density $n_0$. The relation above is obtained by realizing that the function $t \mapsto j[e^{-i t H} \vartheta[n] e^{i tH}]$ is odd when the current of $\vartheta[n]$ is zero. As a consequence, after expanding the exponentials, it follows that the commutators that are even w.r.t. $t$ must vanish. We will use the notation $a \lesssim b$ when there exists $C>0$ such that $a \leq C b$. 

By a Taylor-Lagrange expansion, we observe that, for any $r\in\mathbb{R}$ and $\sigma\in\mathcal{E}^2$,
\begin{equation}\label{eq:expBraH2}
e^{-i \varepsilon r H}\sigma e^{i\varepsilon r H} = \sigma - i\varepsilon \int_0^{r}e^{-i \varepsilon u H}[H,\sigma] e^{i \varepsilon u H} du.
\end{equation}
The following bounds on the current will be utilized throughout the section: denoting by $(\lambda_p,\psi_p)_{p \in \NN}$ the spectral elements of $\sigma$, we have
\begin{equation*}
\nabla\cdot j[\sigma] = 2\Im\left(\sum_{p\in\NN}\lambda_p \psi_p^*(x)\Delta \psi_p(x)  \right),
\end{equation*}
and, for any $t\in\mathbb{R}$,
\bee
\nabla\cdot j[e^{i t H}\sigma e^{-i t H}] &=& 2\Im\left(\sum_{p\in\NN}\lambda_p e^{-i t H}\psi_p^*(x)\Delta e^{i t  H}\psi_p(x)  \right)\\
&=&2\Im\left(\sum_{p\in\NN}\lambda_p e^{-i t H}\psi_p^*(x)e^{i t  H}(\Delta -V)\psi_p(x)  \right).
\eee
Using the Cauchy-Schwarz inequality and the fact that $e^{-i t H}$ is an isometry on $L^2(\Omega)$, we find 
\begin{equation}\label{eq:estimj}
\|\nabla\cdot j[e^{-i t H}\sigma e^{it H}]\|_{L^1} \leq 2\|\sigma\|_{\calJ_1}^{1/2}\|\sigma\|_{\mathcal{E}^2}^{1/2} \leq 2\|\sigma\|_{\mathcal{E}^2}.
\end{equation}
The following estimate will also be helpful in the analysis: for $\sigma \geq 0$ and $V \in L^\infty(\Omega)$,
\be \label{EstE2}
\Tr \big( H_0 e^{-i t H} \sigma e^{i t H} H_0 \big) \lesssim \| \sigma \|_{\calE^2}.
\ee
It is trivial when $H=H_0$, and is obtained as follows when $H=H_0+V$: since
\bee
H_0 e^{-it H} \sigma e^{i t H} H_0&=&H_0 e^{-i t H} (H_0+I)^{-1}(H_0+I)\sigma (H_0+I) (H_0+I)^{-1}e^{i t H}  H_0,
\eee
it suffices to show that  $(H_0+I)^{-1}e^{i t H}  H_0$ is a bounded operator ($I$ denotes the identity).  We write for this
\bee
(H_0+I)^{-1}e^{i t H}  H_0&=&(H_0+I)^{-1}e^{i t H}  H-(H_0+I)^{-1}e^{i t H}V\\
&=&(H_0+I)^{-1}He^{i t H}-(H_0+I)^{-1}e^{i t H}V\\
&=&(H_0+I)^{-1}H_0e^{i t H}+(H_0+I)^{-1}V e^{i t H}-(H_0+I)^{-1}e^{i tH}V,
\eee
and the conclusion follows directly from the facts that $V \in L^\infty(\Omega)$ and $e^{-i t H}$ is unitary.

As explained at the beginning of this Section, there is no general theory for QLE, QDD, and the map $n_0 \mapsto \vartheta[n_0]$. We therefore make the following  assumptions that allow us to derive the error estimates.

  \begin{assumption}\label{hyp}

The following hypotheses will be in force:
\begin{enumerate}
\item (Well-posedness and stability for the modified QDD equation) For any $\varepsilon\geq0$, $T>0$, $s\geq0$ and $\mathsf{n}_s,\mathsf{f}\in L^1$, the following PDE
\begin{equation*}
\mathsf{n}(t) -  \mathsf{n}_s +  \int_s^t \kappa_{\varepsilon}(r-s)\nabla \cdot \left(\mathsf{n}(r) \nabla \mathsf{A}[\mathsf{n}(r)]\right) dr = \int_s^t \mathsf{f}(r)dr,
\end{equation*}
admits a unique solution $\mathsf{n}\in\mathcal{C}([s,s+T]; L^1)$. We denote by $\mathcal{Q}$ the (nonlinear) propagator associated to this PDE (\textit{i.e.} the solution for $t \geq s$ reads $\mathsf{n}(t) = \mathcal{Q}^{\varepsilon}_{s,t}(\mathsf{n}_s,\mathsf{f})$) and we assume that there exists constants $C$ and $\nu\in\mathbb{R}^+$, independent of $\varepsilon$ and $t>s\geq 0$, such that, for any $\mathsf{n}_{1,s},\mathsf{n}_{2,s}\in L^1$ and $\mathsf{f}_1,\mathsf{f}_2\in L^1([0,T];L^1)$,
\begin{align}\label{eq:hyp1}
\|\mathcal{Q}^{\varepsilon}_{s,t}(\mathsf{n}_{1,s},\mathsf{f}_1) &- \mathcal{Q}^{\varepsilon}_{s,t}(\mathsf{n}_{1,s},\mathsf{f}_2)\|_{L^1} \nonumber
\\ \hspace{2em}&\leq \left(\|\mathsf{n}_{1,s} - \mathsf{n}_{2,s}\|_{L^1} + C \|\mathsf{f}_1 - \mathsf{f}_2\|_{L^1([s,t];L^1)}\right)e^{\nu(t-s)}.
\end{align}
\item (Uniform bounds for the full solution) Let $\mn^\eps(t)$ be the density associated with the solution to QLE \fref{eq:Liouville-BGK}. We suppose there exists a constant $C>0$ independent of $\eps$ such that, for all $t \in [0,T]$,
\begin{align} \nonumber
  \|\partial_t \vartheta[\mn^\eps](t)\|_{\calE^2}+\|[H,\partial_t \vartheta[\mn^\eps](t)]\|_{\calE^2}&+\|[H,[H,\vartheta[\mn^\eps](t)]]]\|_{\calE^2} \\
  &+\|[H,[H,[H,\vartheta[\mn^\eps](t)]]]\|_{\calE^2} \leq C \label{eq:hyp1b}
\end{align}
\item (Stability for the equilibrium) There exists a constant $C>0$ such that
\begin{align} \nonumber
  \|\vartheta[\mathsf{n}_1] - \vartheta[\mathsf{n}_2]\|_{\calE^2}&+\|[H,\vartheta[\mathsf{n}_1] - \vartheta[\mathsf{n}_2]]\|_{\calE^2}\\
  &+\|[H,[H,\vartheta[\mathsf{n}_1] - \vartheta[\mathsf{n}_2]]]\|_{\calE^2} \leq C \|\mathsf{n}_1- \mathsf{n}_2 \|_{L^1}.\label{eq:hyp2}
\end{align}
\end{enumerate}
\end{assumption}

Beyond the assumptions above, we will suppose in addition that the semi-discrete version of the modified QDD equation (i.e. the first equation in \fref{dQDD}), enjoys the same stability and regularity properties as the continuous in time version.

We have all needed now to proceed to the derivation of uniform error estimates, and our main results is the following:

\begin{theorem} \label{th} Fix $T>0$, and let $N \in \NN$ such that $N \delta t \leq T$. We have then the following uniform (in $\eps$) estimates between the solutions to QLE and those of the semi-discrete AP scheme \fref{dQDD}:
  $$
\|\mn^\eps(t_n) - \mn^\eps_n \|_{L^1} +\|\varrho^{\varepsilon}(t_{n})-\varrho^{\varepsilon}_{n} \|_{\calE^2}\lesssim \delta t, \qquad 0 \leq n \leq N.
$$
  \end{theorem}

  The proof of the Theorem begins with the proof of Lemma \ref{lem:expsigma}.

\subsection{Proof of Lemma \ref{lem:expsigma}.} The starting point is to derive the expression of the term $\varsigma^{\varepsilon}_1(s,t)$. For this, we rewrite \eqref{eq:Liouville-BGK} in the mild formulation (we drop the $\eps$ in $\varrho^\eps$ for simplicity)
\begin{align} \nonumber
\varrho(t) &= e^{-(t-s)/\varepsilon^2}e^{-i(t-s)H/\varepsilon}\varrho(s)e^{i(t-s)H/\varepsilon} + \varepsilon^{-2}\int_s^te^{-(t-r)/\varepsilon^2} e^{-i(t-r)H/\varepsilon}\vartheta[\mathsf{n}^{\varepsilon}](r)e^{i(t-r)H/\varepsilon} dr
\\ &= S_{\varepsilon,(t-s)/\varepsilon^2}[\varrho(s)] + \int_0^{\varepsilon^{-2}(t-s)}S_{\varepsilon,\tau}[\vartheta[\mn^{\varepsilon}](t-\varepsilon^{2}\tau)] d\tau, \label{mild}
\end{align}
where $S_{\varepsilon,\tau}$ is defined in \fref{defkap}. 
We will use the following two approximations, which are simply Taylor-Lagrange formulas:
\begin{equation*}
\vartheta[\mathsf{n}^{\varepsilon}](t-\varepsilon^{2}\tau) = \vartheta[\mn^{\varepsilon}](t) - \int_{t-\varepsilon^2\tau}^t \partial_t \vartheta[\mn^{\varepsilon}](r)dr,
\end{equation*}
and, for any density operator $\sigma$,
  $$
e^{-i \varepsilon \tau H}\sigma e^{i\varepsilon \tau H} = \sigma - i\tau \varepsilon [H,\sigma] - \varepsilon^2\int_0^{\tau}e^{-i \varepsilon r H}[H,[H,\sigma]] e^{i \varepsilon r H}(\tau-r) dr.
$$
Using both approximations yields
\begin{align*}
\int_0^{\varepsilon^{-2}(t-s)} & S_{\varepsilon,\tau}[\vartheta[\mn^{\varepsilon}](t-\varepsilon^{2}\tau)]  d\tau 
\\ &= \int_0^{\varepsilon^{-2}(t-s)}e^{-\tau} e^{-i\varepsilon \tau H}\vartheta[\mn^{\varepsilon}](t)e^{i\varepsilon \tau H} d\tau -\int_0^{\varepsilon^{-2}(t-s)}S_{\varepsilon,\tau}\left[\int_{t-\varepsilon^2\tau}^t \partial_t \vartheta[\mn^{\varepsilon}](r)dr\right] d\tau
\\ &=\int_0^{\varepsilon^{-2}(t-s)}e^{-\tau}\left( \vartheta[\mn^{\varepsilon}](t) - i\tau \varepsilon [H,\vartheta[\mn^{\varepsilon}](t) ] \right)d\tau
\\ &\hspace{1em}-\varepsilon^2\int_0^{\varepsilon^{-2}(t-s)}e^{-\tau}\left(\int_0^{\tau}e^{-i \varepsilon r H}[H,[H,\vartheta[\mn^{\varepsilon}](t)]]e^{i \varepsilon r H}(\tau-r)  dr\right)d\tau
\\ &\hspace{1em}-\int_0^{\varepsilon^{-2}(t-s)}S_{\varepsilon,\tau}\left[\int_{t-\varepsilon^2\tau}^t \partial_t \vartheta[\mn^{\varepsilon}](r)dr\right] d\tau.
\end{align*}
Furthermore, a simple calculation shows that
\begin{align*}
\int_0^{\varepsilon^{-2}(t-s)}&e^{-\tau}\left( \vartheta[\mn^{\varepsilon}](t) - i\tau \varepsilon [H,\vartheta[\mn^{\varepsilon}](t) ] \right)d\tau
  \\ &= \vartheta[\mn^{\varepsilon}](t)\left(1-e^{-\varepsilon^{-2}(t-s)}\right) - i\varepsilon [H,\vartheta[\mn^{\varepsilon}](t) ]\left(1 - e^{-\varepsilon^{-2}(t-s)} - \varepsilon^{-2}(t-s) e^{-\varepsilon^{-2}(t-s)}\right)
       \\ &=\vartheta[\mn^{\varepsilon}](t)\left(1-e^{-\varepsilon^{-2}(t-s)}\right) - i\varepsilon [H,\vartheta[\mn^{\varepsilon}](t) ] \kappa_\eps(t-s),
\end{align*}
where $\kappa_\eps$ is defined in \fref{defkap}. With
\begin{align*}
\varsigma^{\varepsilon}_1(s,t) :&= -\varepsilon^2\int_0^{\varepsilon^{-2}(t-s)}e^{-\tau}\left(\int_0^{\tau}e^{-i \varepsilon r H}[H,[H,\vartheta[\mathrm{n}^{\varepsilon}](t)]]e^{i \varepsilon r H}(\tau-r)  dr\right)d\tau
  \\ &\hspace{1em} -\int_0^{\varepsilon^{-2}(t-s)}S_{\varepsilon,\tau}\left[\int_{t-\varepsilon^2\tau}^t \partial_t \vartheta[\mathrm{n}^{\varepsilon}](r)dr\right] d\tau,\\
  :&=\varsigma^{\varepsilon}_{11}(s,t)+\varsigma^{\varepsilon}_{12}(s,t),
\end{align*}
we obtain the following expression of $\varrho \equiv\varrho^\eps$, for $t \geq s$,
$$
\varrho(t)=S_{\varepsilon,(t-s)/\varepsilon^2}[\varrho(s)]+\vartheta[\mn^{\varepsilon}](t)(1-e^{-\varepsilon^{-2}(t-s)}) - i\varepsilon [H,\vartheta[\mn^{\varepsilon}](t) ] \kappa_\eps(t-s)+\varsigma^{\varepsilon}_1(s,t).
$$

\medskip

We now bound $\varsigma^\eps_1$. Consider first the term $\varsigma^{\varepsilon}_{11}$ and decompose $-[H,[H,\vartheta[\mathrm{n}^{\varepsilon}](t)]]$ into positive and negative parts as $-[H,[H,\vartheta[\mathrm{n}^{\varepsilon}](t)]]=\sigma_+-\sigma_-$, where $\sigma_{\pm}\geq 0$ and $|[H,[H,\vartheta[\mathrm{n}^{\varepsilon}](t)]]|=\sigma_++\sigma_-$. The term $\varsigma^{\varepsilon}_{11}$ is then split into $\varsigma^{\varepsilon}_{11+}-\varsigma^{\varepsilon}_{11-}$ where $\varsigma^{\varepsilon}_{11\pm} \geq 0$. Hence,
$$
\Tr \big( H_0 \varsigma^{\varepsilon}_{11\pm} H_0 \big) =\varepsilon^2\int_0^{\varepsilon^{-2}(t-s)}e^{-\tau}\left(\int_0^{\tau}\Tr \big( H_0  e^{-i \varepsilon r H} \sigma_{\mp}e^{i \varepsilon r H}H_0  \big) (\tau-r)  dr\right)d\tau,
$$
and we can use \fref{EstE2} to estimate the trace in the integral. We proceed in the same way with $\varsigma^{\varepsilon}_{12}$ to arrive at 
\bee
\|\varsigma^{\varepsilon}_1(s,t)\|_{\calE^2} &\leq & \eps^2 \|[H,[H,\vartheta[\mathrm{n}^{\varepsilon}](t)]] \|_{\calE^2} \int_0^{\varepsilon^{-2}(t-s)}e^{-\tau} \int_0^{\tau} (\tau-r)  dr d\tau \\
&&+\eps^2 \sup_{r \in \Rm_+ }\|\partial_t \vartheta[\mn^{\varepsilon}](r) \|_{\calE^2}\int_0^{\varepsilon^{-2}(t-s)}\tau e^{-\tau} d \tau.
\eee
Using the same method, we find as well
\bee
\|[H,\varsigma^{\varepsilon}_1(s,t)]\|_{\calE^2} &\leq & \eps^2 \|[H,[H,[H,\vartheta[\mathrm{n}^{\varepsilon}](t)]]] \|_{\calE^2} \int_0^{\varepsilon^{-2}(t-s)}e^{-\tau} \int_0^{\tau} (\tau-r)  dr d\tau \\
&&+\eps^2 \sup_{r \in \Rm_+ }\| [H,\partial_t \vartheta[\mn^{\varepsilon}](r)] \|_{\calE^2}\int_0^{\varepsilon^{-2}(t-s)}\tau e^{-\tau} d \tau.
\eee
The integrals w.r.t. $\tau$ above can be bounded by $C (\eps^{-4} (t-s)^2) \wedge 1$, so that
\begin{align*}
  \|[H,\varsigma^{\varepsilon}_1(s,t)]\|_{\calE^2} \lesssim&\left( \sup_{r \in \Rm_+ }\|[H,[H,[H,\vartheta[\mathrm{n}^{\varepsilon}](r )]]] \|_{\calE^2} +\sup_{r \in \Rm_+ }\| [H,\partial_t \vartheta[\mn^{\varepsilon}](r)] \|_{\calE^2} \right)\\
  & \hspace{2cm}\times \eps^2 \wedge ( \eps^{-2}(t-s)^2),
\end{align*}
with a similar estimate for $\|\varsigma^{\varepsilon}_1(s,t)\|_{\calE^2}$. 
Using \fref{eq:hyp1b}, it follows that
\be \label{estvarS1}
 \|\varsigma^{\varepsilon}_1(s,t)\|_{\calE^2}+\|[H,\varsigma^{\varepsilon}_1(s,t)]\|_{\calE^2} \lesssim \eps^2 \wedge ( \eps^{-2}(t-s)^2).
\ee
This estimate is sufficient to handle the density operator $\varrho^\eps$, but, because of the additional $\eps^{-1}$, is not sufficient for the current $\eps^{-1} j[\varsigma^{\varepsilon}_1(s,t)]$. We need then to refine the expression of $\varsigma^{\varepsilon}_1$ as follows: we use \eqref{eq:expBraH2} to arrive at
\begin{align*}
  \varsigma^{\varepsilon}_1(s,t) &= -\varepsilon^2\int_0^{\varepsilon^{-2}(t-s)}e^{-\tau}\left(\int_0^{\tau}[H,[H,\vartheta[\mathrm{n}^{\varepsilon}](t)]](\tau-r)  dr\right)d\tau\\
                                 &\hspace{1em}+i\varepsilon^3\int_0^{\varepsilon^{-2}(t-s)}e^{-\tau}\left(\int_0^{\tau}\int_0^{r}e^{-i \varepsilon u H}[H,[H,[H,\vartheta[\mathrm{n}^{\varepsilon}](t)]]]e^{i \varepsilon u H}(\tau-r)dudr \right) d\tau
                                    \\ &\hspace{1em} -\int_0^{\varepsilon^{-2}(t-s)}e^{-\tau}\left[\int_{t-\varepsilon^2\tau}^t \partial_t \vartheta[\mathrm{n}^{\varepsilon}](r)dr\right] d\tau
\\ &\hspace{1em}  +i\eps \int_0^{\varepsilon^{-2}(t-s)}e^{-\tau}\int_0^\tau e^{-i \varepsilon u H} \left[\int_{t-\varepsilon^2\tau}^t [H,\partial_t \vartheta[\mathrm{n}^{\varepsilon}](r)]dr\right] e^{i \varepsilon u H} du d\tau
\end{align*}

The key observation is that the terms of order $\eps^2$ above, namely the first and third terms, have zero current because of \fref{commu} and the fact that $j[\vartheta[\mathrm{n}^{\varepsilon}](t)]=0$. Hence,
\begin{align*}
  j[\varsigma^{\varepsilon}_1(s,t)] &= i\varepsilon^3\int_0^{\varepsilon^{-2}(t-s)}e^{-\tau}\left(\int_0^{\tau}\int_0^{r}j\left[e^{-i \varepsilon u H}[H,[H,[H,\vartheta[\mathrm{n}^{\varepsilon}](t)]]]e^{i \varepsilon u H}\right](\tau-r)dudr \right) d\tau
\\ &\hspace{1em}  +i\eps \int_0^{\varepsilon^{-2}(t-s)}e^{-\tau}\int_0^\tau j\left[e^{-i \varepsilon u H} \left[\int_{t-\varepsilon^2\tau}^t [H,\partial_t \vartheta[\mathrm{n}^{\varepsilon}](r)]dr\right] e^{i \varepsilon u H} \right]du d\tau.
\end{align*}
Denoting by $j_0(t,u)$ the first current in the r.h.s. above, and by $j_1(r,u)$ the term $j\left[e^{-i \varepsilon u H} [H,\partial_t \vartheta[\mathrm{n}^{\varepsilon}](r)]e^{i \varepsilon u H} \right]$, we find the estimate
$$
\|\nabla \cdot j[\varsigma^{\varepsilon}_1(s,t)]\|_{L^1} \leq C \eps \left(\sup_{u \in \Rm_+ } \|\nabla \cdot j_0(t,u) \|_{L^1}+\sup_{r,u \in \Rm_+ } \|\nabla \cdot j_1(r,u) \|_{L^1} \right) \eps^2 \wedge ( \eps^{-2}(t-s)^2).
$$
Using \fref{eq:hyp1b} and \fref{eq:estimj} to bound the currents $j_0$ and $j_1$, we find
\be \label{estvarS2}
 \|\nabla \cdot j[\varsigma^{\varepsilon}_1(s,t)]\|_{L^1} \lesssim \eps [\eps^2 \wedge ( \eps^{-2}(t-s)^2)].
\ee
Compared to \fref{estvarS1}, we have gained a factor $\eps$ for the current of $\varsigma^{\varepsilon}_1$ (more precisely its divergence). Observing that $\eps^2 \wedge ( \eps^{-2}(t-s)^2)=O(\eps^2 \wedge (t-s))$, this ends the proof of Lemma \ref{lem:expsigma}.

\subsection{Proof of Theorem \ref{th}.} The proof is divided into two steps. In the first one, we neglect the term $\varsigma_1^\eps$ in QLE  to arrive at the intermediate scheme \fref{eq:scheme1}. We will see that this introduces an error of order $O (\eps^2 \wedge \delta t)$. In the second step, we discretize the time integrals and approximate the propagator $e^{-i tH}$at second order globally. This will generate an error of order $O(\delta t)$, leading to an overall error $O(\delta t)$ that is therefore uniform in $\eps$.

\paragraph{Step 1: neglecting $\varsigma_1^\eps$.} We recall that $(\tilde \mn^\eps,\tilde \varrho_n^\eps)$ denote the solution to the system \fref{eq:scheme1} where $\varsigma_1^\eps$ is removed.
Our goal in this Section is to estimate the errors 
\begin{equation*}
\boldsymbol{\eta}^{\varepsilon}_n : =\tilde \mn^{\varepsilon}(t_n) -  n[\varrho^{\varepsilon}(t_n)], \qquad \boldsymbol{\theta}_n^{\varepsilon} : = \tilde \varrho^{\varepsilon}_n - \varrho^{\varepsilon}(t_n),
\end{equation*}
where $n[\varrho^{\varepsilon}(t_n)]$ and $\varrho^{\varepsilon}(t_n)$ are the exact solutions at time $t_n$. Despite its largely informal character, the proof is quite technical and is based on controlling  $\varsigma_1^\eps$ and remainders, and on the stability estimates in Assumptions \ref{hyp}.

According to Lemma \ref{lem:expsigma}, we have the expression
\begin{equation*}
\boldsymbol{\theta}_{n+1}^{\varepsilon} = S_{\varepsilon,\delta t/\varepsilon^2}[\boldsymbol{\theta}_n^{\varepsilon}] +  \boldsymbol{\nu}_{n+1}^{\varepsilon}\left(1-e^{-\varepsilon^{-2}\delta t}\right) - i\varepsilon \kappa_{\varepsilon}(\delta t) [H,\boldsymbol{\nu}_{n+1}^{\varepsilon} ] - \varsigma^{\varepsilon}_1(t_n,t_{n+1}),
\end{equation*}
where 
\begin{equation*}
\boldsymbol{\nu}_n^{\varepsilon} = \vartheta[\tilde \mn^{\varepsilon}(t_{n})] - \vartheta[n[\varrho^{\varepsilon}(t_n)]].
\end{equation*}
Since $\boldsymbol{\theta}_{0}^{\varepsilon} = 0$, we find, by iterating over $\boldsymbol{\theta}_{n}^{\varepsilon}$, for any $n\geq 1$,
\begin{equation} \label{eqthet}
\boldsymbol{\theta}_{n}^{\varepsilon} = \sum_{k = 1}^n S_{\varepsilon, (n-k)\delta t/\varepsilon^2}\left[ \boldsymbol{\nu}_{k}^{\varepsilon}\left(1-e^{-\varepsilon^{-2}\delta t}\right) - i\varepsilon \kappa_{\varepsilon}(\delta t) [H,\boldsymbol{\nu}_{k}^{\varepsilon} ] - \varsigma^{\varepsilon}_1(t_{k-1},t_{k}) \right],
\end{equation}
so that, for any $r\in [t_n,t_{n+1}]$,
\begin{equation} \label{eqthet2}
S_{\varepsilon,(r-t_n)/\varepsilon^2}[\boldsymbol{\theta}_{n}^{\varepsilon}] = \sum_{k = 1}^n S_{\varepsilon, (r- k \delta t)/\varepsilon^2}\left[ \boldsymbol{\nu}_{k}^{\varepsilon}\left(1-e^{-\varepsilon^{-2}\delta t}\right) - i\varepsilon \kappa_{\varepsilon}(\delta t) [H,\boldsymbol{\nu}_{k}^{\varepsilon} ] - \varsigma^{\varepsilon}_1(t_{k-1},t_{k}) \right].
\end{equation}
According to the first equation in \fref{eq:scheme1}, we need to estimate currents in the r.h.s. For this, from \fref{eq:estimj} and \fref{eq:expBraH2}, we will exploit that
\be \label{estJ2}
\| \nabla \cdot j[e^{-i  r H}\sigma e^{i r H}] \|_{L^1} \leq \| \nabla \cdot j[\sigma]\|_{L^1} + C r \| [H,\sigma] \|_{\calE^2}
\ee
and also
$$
\| \nabla \cdot j[e^{-i  r H}\sigma e^{i r H}] \|_{L^1} \leq \|\sigma \|_{\calE^2}.
$$
We will use below crucially the facts that $j[ \boldsymbol{\nu}_{k}^{\varepsilon}]=0$ and, as seen in the proof of Lemma \ref{lem:expsigma}, that $j[\varsigma^{\varepsilon}_1]$ behaves better than $\varsigma^{\varepsilon}_1$ w.r.t. $\eps$. This is exploited in combination with \fref{estJ2}, and we find the estimates
\begin{align*}
&\| \nabla \cdot j[S_{\varepsilon, (r- k \delta t)/\varepsilon^2}[ \boldsymbol{\nu}_{k}^{\varepsilon}] \|_{L^1} \lesssim  \frac{(r-k \delta t)}{\eps} e^{- (r-k \delta t)/\eps^2} \| [H,\boldsymbol{\nu}_{k}^{\varepsilon}] \|_{\calE^2}\\
  &\| \nabla \cdot j[S_{\varepsilon,(r- k \delta t)/\varepsilon^2}[\varsigma^{\varepsilon}_1(t_{k-1},t_{k}) ] \|_{L^1} \lesssim e^{-\frac{r-k\delta t}{\varepsilon^2}} \|  \nabla \cdot j[\varsigma^{\varepsilon}_1(t_{k-1},t_{k})]\|_{L^1}\\
  &\hspace{7cm}+\frac{(r-k \delta t)}{\eps} e^{- (r-k \delta t)/\eps^2}\| [H,\varsigma^{\varepsilon}_1(t_{k-1},t_{k})] \|_{\calE^2}.
\end{align*}
Going back to \fref{eqthet2} and utilizing these last two estimates, it follows that
\begin{align*}
  \|\nabla \cdot& j[S_{\varepsilon,(r-t_n)/\varepsilon^2}[\boldsymbol{\theta}_{n}^{\varepsilon}]]\|_{L^1}\\
  \lesssim &  \sum_{k = 1}^n e^{- (r-k \delta t)/\eps^2}\left[  \left(\frac{(r-k \delta t)}{\eps} \left(1-e^{-\varepsilon^{-2}\delta t}\right) + \eps \kappa_\eps(\delta t)\right)\| [H,\boldsymbol{\nu}_{k}^{\varepsilon}] \|_{\calE^2}\right]\\
&+ \sum_{k = 1}^n e^{- (r-k \delta t)/\eps^2} \left( \|  \nabla \cdot j[\varsigma^{\varepsilon}_1(t_{k-1},t_{k})]\|_{L^1} + \frac{r-k\delta t}{\varepsilon} \| [H,\varsigma^{\varepsilon}_1(t_{k-1},t_{k})] \|_{\calE^2} \right).
\end{align*}
Since $\kappa_\eps \leq 1-e^{-\varepsilon^{-2}\delta t}$ and $\frac{(r-k \delta t)}{\eps^2} e^{-\frac{r-k\delta t}{\varepsilon^2}} \lesssim  e^{-\frac{r-k\delta t}{2\varepsilon^2}}$, this reduces to
\begin{align} \nonumber
  \eps^{-1}\|\nabla \cdot &j[S_{\varepsilon,(r-t_n)/\varepsilon^2}[\boldsymbol{\theta}_{n}^{\varepsilon}]]\|_{L^1}\lesssim \sum_{k = 1}^n e^{- (r-k \delta t)/2\eps^2} \left(1-e^{-\varepsilon^{-2}\delta t}\right)\| [H,\boldsymbol{\nu}_{k}^{\varepsilon}] \|_{\calE^2}\\
&+ \sum_{k = 1}^n e^{- (r-k \delta t)/2\eps^2} \left( \eps^{-1}\|  \nabla \cdot j[\varsigma^{\varepsilon}_1(t_{k-1},t_{k})]\|_{L^1} + \| [H,\varsigma^{\varepsilon}_1(t_{k-1},t_{k})] \|_{\calE^2} \right).\label{boundthet}
\end{align}
Using Assumption \fref{eq:hyp2} to control $\boldsymbol{\nu}_{k}^{\varepsilon}$ in terms of $\boldsymbol{\eta}_{k}^{\varepsilon}$, we finally obtain the estimate
\begin{align}
  \eps^{-1}\|\nabla \cdot j[&S_{\varepsilon,(r-t_n)/\varepsilon^2}[\boldsymbol{\theta}_{n}^{\varepsilon}]]\|_{L^1}\lesssim 
   \left(1-e^{-\varepsilon^{-2}\delta t}\right) \sum_{k = 1}^n e^{- (r-k \delta t)/2\eps^2} \| \boldsymbol{\eta}_{k}^{\varepsilon} \|_{L^1} \nonumber
\\
&+ \sum_{k = 1}^n e^{- (r-k \delta t)/2\eps^2} \left( \eps^{-1}\|  \nabla \cdot j[\varsigma^{\varepsilon}_1(t_{k-1},t_{k})]\|_{L^1} + \| [H,\varsigma^{\varepsilon}_1(t_{k-1},t_{k})] \|_{\calE^2} \right). \label{estimat0}
\end{align}


We are now ready to use the modified QDD equation and Assumption \fref{eq:hyp2} to estimate $\boldsymbol{\eta}_{n}^{\varepsilon}$. For this, let first
\bee
J_n^\eps&=&\delta t \eps^{-1} \sup_{r} \|\nabla \cdot j[\varsigma^{\varepsilon}_1(t_n,r)]\|_{L^1}\\
&&+\delta t \sum_{k = 1}^n e^{- (n-k) \delta t/2\eps^2} \left( \eps^{-1}\|  \nabla \cdot j[\varsigma^{\varepsilon}_1(t_{k-1},t_{k})]\|_{L^1} + \| [H,\varsigma^{\varepsilon}_1(t_{k-1},t_{k})] \|_{\calE^2} \right).
\eee
With the $\calQ$ notation of Assumptions \ref{hyp}, we have $\mn^\eps(t_{n+1})=\calQ^\eps_{t_n,t}(\mn^\eps(t_{n}),f_1)$, where
$$f_1(t)= \eps^{-1} \nabla \cdot j[S_{\varepsilon,(t-t_n)/\varepsilon^2}[\varrho^\eps(t_n)]] +\eps^{-1} \nabla \cdot j[\varsigma_1^{\eps}(t_{n},t)],$$ as well as
$\tilde \mn^\eps(t_{n+1})=\calQ^\eps_{t_n,t}(\tilde \mn^\eps(t_{n}),f_2)$, where $f_2(t)= \eps^{-1} \nabla \cdot j[S_{\varepsilon,(t-t_n)/\varepsilon^2}[\tilde \varrho^\eps(t_n)]]$. Since $e^{- (r-k \delta t)/2\eps^2}\leq e^{- (n-k) \delta t/2\eps^2}$ when $r \in [t_n,t_{n+1}]$, we have, combining \fref{estimat0} with Assumptions \fref{eq:hyp1} and \fref{boundthet},
\begin{align*}
\|\boldsymbol{\eta}^{\varepsilon}_{n+1} \|_{L^1} &\leq \left(\|\boldsymbol{\eta}^{\varepsilon}_n\|_{L^1} + \varepsilon^{-1}C\int_{t_n}^{t_{n+1}}\|\nabla \cdot j[S_{\varepsilon,(r-t_n)/\varepsilon^2}[\boldsymbol{\theta}_{n}^{\varepsilon}]]\|_{L^1} dr \right.
\\ &\hspace{1em}\left.+ \varepsilon^{-1}C\int_{t_n}^{t_{n+1}} \|\nabla \cdot j[\varsigma^{\varepsilon}_1(t_n,r)]\|_{L^1}dr\right)e^{\nu \delta t}
\\ &\leq e^{\nu\delta t}\|\boldsymbol{\eta}^{\varepsilon}_n\|_{L^1} + C e^{\nu \delta t}\delta t \left(1-e^{-\varepsilon^{-2}\delta t}\right) \sum_{k = 1}^n  e^{-\frac{(n-k)\delta t}{2\varepsilon^2}} \|\boldsymbol{\eta}^{\varepsilon}_{k}\|_{L^1}+e^{\nu\delta t} J_n^\eps.
\end{align*}
The inequality above has the form
$$
\gamma_{n+1}\leq e^{a} \gamma_n+ b \sum_{k = 1}^n  e^{-\frac{(n-k)\delta t}{2\varepsilon^2}} \gamma_k+c_n,
$$
where
\bee
\gamma_n &=&\|\boldsymbol{\eta}^{\varepsilon}_n\|_{L^1}, \qquad a =\nu \delta t\\
c_n&=&e^{\nu\delta t} J_n^\eps, \qquad b=C e^{\nu \delta t}\delta t \left(1-e^{-\varepsilon^{-2}\delta t}\right).
\eee
Iterating, we find for $N \geq 2$, 
$$
\gamma_N \leq e^{(N-1)a}\gamma_1+b \sum_{n=1}^{N-1} e^{a (n-1)} \sum_{k = 1}^n  e^{-\frac{(n-k)\delta t}{2\varepsilon^2}} \gamma_k+\sum_{n=1}^{N-1} e^{a (n-1)}c_n.
$$
Since
$$
\sum_{n=1}^{N-1} \sum_{k = 1}^n  e^{-\frac{(n-k)\delta t}{2\varepsilon^2}} \gamma_k=
\sum_{k=1}^{N-1} \gamma_k \left(\sum_{n = k}^{N-1}  e^{-\frac{(n-k)\delta t}{2\varepsilon^2}} \right) \leq \sum_{k=1}^{N-1} \gamma_k \sum_{n = 0}^{N-1}  e^{-\frac{ n\delta t}{2\varepsilon^2}}= \frac{1-e^{-\frac{N \delta t}{2\varepsilon^2}}}{1-e^{-\frac{\delta t}{2\varepsilon^2}}} \sum_{k=1}^{N-1} \gamma_k,
$$
we have, with $\alpha_N=b e^{a(N-2)}(1-e^{-\frac{N \delta t}{2\varepsilon^2}})/(1-e^{-\frac{\delta t}{2\varepsilon^2}})$,
$$
\gamma_N \leq e^{(N-1)a}\gamma_1+ \alpha_N \sum_{k = 1}^{N-1} \gamma_k+\sum_{n=1}^{N-1} e^{a (n-1)}c_n.
$$
The discrete Gronwall inequality yields finally
$$
\gamma_N\leq e^{(N-1)\alpha_N} \left( e^{(N-1)a}\gamma_1+\sum_{n=1}^{N-1} e^{a (n-1)}c_n\right).
$$
If the maximal simulation time is $T$, $N$ is set to be the largest integer such that $ N \delta t\leq T$. In our case, $b=C e^{\nu \delta t}\delta t (1-e^{-\varepsilon^{-2}\delta t})$, so that $(N-1)\alpha_N \leq C e^{2\nu T} \leq C$, uniformly in $\delta t$, $N$, and $\eps$. Also $c_n=e^{\nu \delta t} J_n^\eps$, and thanks to \fref{estvarS1}-\fref{estvarS2},
$$
|J_n^\eps| \lesssim \delta t \; [\eps^2 \wedge ( \eps^{-2}\delta t^2) ]\left(1+\sum_{k=1}^{n} e^{-\frac{(n-k)\delta t}{2\varepsilon^2}} \right).
$$
Then,
$$
\sum_{n=1}^{N-1} e^{a (n-1)}c_n \leq e^{\nu T} \sum_{n=1}^{N-1} J_n^\eps \lesssim \eps^2 \wedge ( \eps^{-2}\delta t^2) \left(1+(1-e^{-\frac{\delta t}{2\varepsilon^2}})^{-1}\right) \lesssim \eps^2 \wedge \delta t,
$$
where the last inequality is uniform in $\delta t$ and $\eps$.

It remains to treat the error after one iteration, which is controlled by
\bee
\|\boldsymbol{\eta}^{\varepsilon}_{1} \|_{L^1} &\lesssim& \delta t \eps^{-1} \sup_{r} \|\nabla \cdot j[\varsigma^{\varepsilon}_1(t_0,r)]\|_{L^1}+\delta t  \left( \eps^{-1}\|  \nabla \cdot j[\varsigma^{\varepsilon}_1(t_{0},t_{1})]\|_{L^1} + \| [H,\varsigma^{\varepsilon}_1(t_{0},t_{1})] \|_{\calE^2} \right)\\
&\lesssim& \delta t \; \eps^2 \wedge ( \eps^{-2}\delta t^2) \lesssim \delta t \; \eps^2 \wedge \delta t.
\eee
Collecting previous estimates, we finally find that, for all $0 \leq n \leq N$,
\be \label{interm}
\|\boldsymbol{\eta}^{\varepsilon}_{n} \|_{L^1} \lesssim \eps^2 \wedge \delta t.
\ee
Combining the latter, Assumptions \fref{eq:hyp2} and \fref{eqthet}, we obtain
\be \label{intermthet}
\|\boldsymbol{\theta}^{\varepsilon}_{n} \|_{\calE^2} \lesssim \eps^2 \wedge \delta t.
\ee

We have therefore obtained uniform error estimates between the exact system and the one where $\varsigma_1^\eps$ is neglected. The next step is to compare the latter system with the semi-discrete scheme where time integrals and propagators are discretized.

\paragraph{Step 2: error for the semi-discrete scheme.} We recall that the semi-discrete scheme reads
\begin{equation}
\left\{\begin{array}{ll}
\ds & \ds  \mn^{\varepsilon}_{n+1} -   \mn^{\varepsilon}_n+ \left( \int_{t_n}^{t_{n+1}} \kappa_{\varepsilon}(r-t_n)dr \right) \nabla \cdot \left( \mn^{\varepsilon}_n \nabla \mathsf{A}[\mn^\eps_{n+1}]\right)
\ds \\ &\hspace{8em}\ds =  \eps (1-e^{-\frac{\delta t}{\eps^2}}) \nabla \cdot j[U^\eps_{a_\eps} \varrho_n^\eps (U^\eps_{a_\eps})^*] 
\ds \\ [3mm] & \varrho^{\varepsilon}_{n+1} = \hat{S}_{\varepsilon,\delta t/\varepsilon^2}[ \varrho^{\varepsilon}_n] +  \vartheta[ \mn^{\varepsilon}_{n+1}](1-e^{-\varepsilon^{-2}\delta t}) - i\varepsilon \kappa_{\varepsilon}(\delta t) [H,\vartheta[ \mn^{\varepsilon}_{n+1}] ],
\end{array}\right. \label{fullDS}
\end{equation}
where $\hat S_{\varepsilon,\delta t/\varepsilon^2} [\sigma] =e^{-\delta t/\eps^2} U^\eps_{\delta t} \sigma (U^\eps_{\delta t})^*$ for $U^\eps_{\delta t}$ a globally second order approximation of $e^{-i \delta t H /\eps}$. 
The analysis is quite similar to that of Step 1, so we only underline the main differences and do not detail all the calculations. We will repeatedly use the discrete version of Assumptions \ref{hyp} for the discrete modified QDD equation in order to control how the errors propagate over iterations on $n$.

To control the error between $\tilde \mn^{\varepsilon}(t_{n+1})$ and $\mn^{\varepsilon}_{n+1}$, we take the difference between \fref{fullDS} and \fref{eq:scheme1}.  There are two different types  of errors: one type due to the discretization of the time integrals, and one type due to the discretization of $e^{-i t H /\eps}$. We have more precisely
$$
  \begin{array}{ll}
         \ds & \ds \mn^{\varepsilon}_{n+1} -   \mn^{\varepsilon}_n+ \left( \int_{t_n}^{t_{n+1}} \kappa_{\varepsilon}(r-t_n)dr \right) \nabla \cdot \left( \mn^{\varepsilon}_n \nabla \mathsf{A}[\mn^\eps_{n+1}]\right)  \\
 \ds & \ds  -\left(\tilde \mn^{\varepsilon}_{n+1} -   \tilde \mn^{\varepsilon}_n+ \left( \int_{t_n}^{t_{n+1}} \kappa_{\varepsilon}(r-t_n)dr \right) \nabla \cdot \left( \tilde \mn^{\varepsilon}_n \nabla \mathsf{A}[\tilde \mn^\eps_{n+1}]\right)\right)          
\ds \\ &\hspace{8em}\ds =  R_1+R_2+R_3+R_4
\end{array}
$$
where
\begin{align*}
  &R_1= \eps^{-1}\int_{t_n}^{t_{n+1}}e^{-\frac{r-t_n}{\varepsilon^2}} \nabla \cdot j[e^{-i (r-t_n) H /\eps}\tilde \varrho_n^\eps e^{i (r-t_n) H /\eps}]]dr \\
  &\hspace{3cm} -  \eps^{-1}\int_{t_n}^{t_{n+1}}e^{-\frac{r-t_n}{\varepsilon^2}} \nabla \cdot j[e^{-i a_\eps H /\eps}\tilde \varrho_n^\eps e^{ia_\eps H /\eps}]]dr\\
  &R_2=\eps^{-1}\int_{t_n}^{t_{n+1}}e^{-\frac{r-t_n}{\varepsilon^2}} \nabla \cdot j[e^{-i a_\eps H /\eps}\tilde \varrho_n^\eps e^{ia_\eps H /\eps}]]dr\\&\hspace{3cm} - \eps^{-1}\int_{t_n}^{t_{n+1}}e^{-\frac{r-t_n}{\varepsilon^2}} \nabla \cdot j[ e^{-i a_\eps H /\eps}\varrho_n^\eps e^{i a_\eps H /\eps}]dr\\
  &R_3=\eps^{-1}\int_{t_n}^{t_{n+1}}e^{-\frac{r-t_n}{\varepsilon^2}} \nabla \cdot j[e^{-i a_\eps H /\eps}\varrho_n^\eps e^{ia_\eps H /\eps}]]dr\\&\hspace{3cm} - \eps^{-1}\int_{t_n}^{t_{n+1}}e^{-\frac{r-t_n}{\varepsilon^2}} \nabla \cdot j[ U^\eps_{a_\eps}\varrho_n^\eps U^\eps_{a_\eps}]dr\\
  &R_4=\int_{t_n}^{t_{n+1}} \kappa_{\varepsilon}(r-t_n) \nabla \cdot \left( \tilde \mn^{\varepsilon}(r) \nabla \mathsf{A}[\tilde \mn^\eps(r)] \right)dr\\
  &\hspace{3cm}-\int_{t_n}^{t_{n+1}} \kappa_{\varepsilon}(r-t_n)dr \nabla \cdot \left(\tilde \mn^{\varepsilon}_n \nabla \mathsf{A}[\tilde \mn^\eps_{n+1}]\right).
  \end{align*}

The term $R_1$ has a similar role as $J_n^\eps$ in the previous Section. According to the Taylor-Lagrange formula, $R_1$ can be expressed as
\be \label{2nd}
\eps^{-1}\int_{t_n}^{t_{n+1}}dr e^{-\frac{r-t_n}{\varepsilon^2}} \int_{a_\eps/\eps}^{(r-t_n)/\eps }\nabla \cdot j\left[e^{-i \tau H}[H,[H,\tilde \varrho_n^\eps]] e^{i \tau H}\right]((r-t_n)/\eps-\tau) d \tau.
\ee
Note that the first order term vanishes because of our choice for $a_\eps$. To estimate the term above, we need to control the current. For this, iterating the second equation in \fref{eq:scheme1}, we find that the operator $\tilde \varrho_{n}^{\varepsilon}$ verifies
\begin{equation*}
\tilde \varrho_{n}^{\varepsilon} = S_{\varepsilon, n \delta t/\varepsilon^2} [\varrho_0]+\sum_{k = 1}^n S_{\varepsilon, (n-k)\delta t/\varepsilon^2}\left[ \vartheta[\tilde \mn^{\varepsilon}](t_k)\left(1-e^{-\varepsilon^{-2}\delta t}\right) - i\varepsilon \kappa_{\varepsilon}(\delta t) [H,\vartheta[\tilde \mn^{\varepsilon}](t_k)]\right],
\end{equation*}
which accordingly allows us to decompose $j\left[e^{-i \tau H}[H,[H,\tilde \varrho_n^\eps]] e^{i \tau H}\right]$ into 3 terms, one linear in $\varrho_0$, one linear in $\vartheta[\tilde \mn^{\varepsilon}](t_k)$, and one linear in the commutator $[H,\vartheta[\tilde \mn^{\varepsilon}](t_k)]$. After integration over $\tau$ as in \fref{2nd} and taking the $L^1$ norm, the first and last term in this decomposition are bounded, thanks to \fref{eq:estimj}, by
$$
\frac{(r-t_n)^2}{\eps^2} \left(e^{- n \delta t/\eps^2}\| [H,[H,\varrho_{0}]]\|_{\calE^2}+ \eps \sup_{k\leq n} \| [H,[H,[H,\vartheta[\tilde \mn^\eps](t_k)]]\|_{\calE^2}\right).
$$
For the second term of the decomposition, we crucially use once more that
$$j[[H,[H,\vartheta[\tilde \mn^\eps](t_k)]]]=0,$$ and together with the fact that $[H,S_{\varepsilon,\tau}[\sigma]] = S_{\varepsilon,\tau}[[H,\sigma]]$, for any $\tau\geq 0$ and $\sigma\in\mathcal{E}^2$, as well as \fref{eq:expBraH2}, we find a bound of the form
$$
\frac{(r-t_n)^3}{\eps^3} \sup_{k \leq n} \| [H,[H,[H,\vartheta[\tilde \mn^\eps](t_k)]]]\|_{\calE^2}.
$$
Assuming as in \fref{eq:hyp1b} that the various $\|\cdot \|_{\calE^2}$ norms above are uniformly bounded in $\eps$ and $n$, we find that
$$
\|R_1\|_{L^1} \lesssim e^{- n \delta t/\eps^2} \int_{t_n}^{t_{n+1}}e^{-\frac{r-t_n}{\varepsilon^2}} \frac{(r-t_n)^2}{\eps^3}  dr+\int_{t_n}^{t_{n+1}}e^{-\frac{r-t_n}{\varepsilon^2}} \frac{(r-t_n)^2}{\eps^2} \left( 1+\frac{(r-t_n)}{\eps^2} \right)dr.
$$
After summing over $n$ up to $n=N$, the second term above is an $O(\eps^2 \wedge \delta t)$, while the first one is controlled by
$$
(1-e^{-\delta t/\eps^2})^{-1}\int_{0}^{\delta t }e^{-\frac{r}{\varepsilon^2}} \frac{r^2}{\eps^3}  dr \lesssim \eps (\eps^2 \wedge \delta t).
$$
Above, we used that $(1-e^{-\delta t/\eps^2})^{-1} \delta t /\eps^2 \lesssim 1$. The conclusion is that the sum of $R_1$ over $n$ has the same order as the term $\sum_{n=1}^{N-1} J_n^\eps$ of the previous Section, which is an $O(\eps^2 \wedge \delta t)$. The error introduced by $R_1$ on $\mn^{\varepsilon}_{n}-\tilde \mn^\eps$ is hence of order  $O( \eps^2 \wedge \delta t)$.

The term $R_2$ is treated in a similar manner as $\eps^{-1} \nabla \cdot j[S_{\varepsilon,(r-t_n)/\varepsilon^2}[\boldsymbol{\theta}_{n}^{\varepsilon}]]$ in the previous Section. Indeed, with
\begin{equation*}
 \boldsymbol{\tilde \eta}^{\varepsilon}_n : = \tilde \mn^{\varepsilon}_n -  \mn^{\varepsilon}_{n}, \qquad \boldsymbol{\tilde \theta}_n^{\varepsilon} : =  \tilde \varrho^{\varepsilon}_n - \varrho^{\varepsilon}_n,
\end{equation*}
we have
\begin{equation*}
\boldsymbol{\tilde \theta}_{n+1}^{\varepsilon} = S_{\varepsilon,\delta t/\varepsilon^2}[\boldsymbol{\tilde \theta}_n^{\varepsilon}] +  \boldsymbol{\tilde \nu}_{n+1}^{\varepsilon}\left(1-e^{-\varepsilon^{-2}\delta t}\right) - i\varepsilon \kappa_{\varepsilon}(\delta t) [H,\boldsymbol{\tilde \nu}_{n+1}^{\varepsilon} ] - \tilde \varsigma^{\varepsilon}_{1,n},
\end{equation*}
with $\boldsymbol{\tilde \theta}_{0}^{\varepsilon} = 0$ and where
$$
\tilde \varsigma^{\varepsilon}_{1,n}=S_{\varepsilon,\delta t/\varepsilon^2}[\varrho^{\varepsilon}_n]-\hat S_{\varepsilon,\delta t/\varepsilon^2}[\varrho^{\varepsilon}_n],
$$
and
\begin{equation*}
\boldsymbol{\tilde \nu}_n^{\varepsilon} = \vartheta[\tilde \mn^{\varepsilon}_n] - \vartheta[\mn^{\varepsilon}_n].
\end{equation*}
Since $U^\eps_{t}$ is a globally second order approximation, we have
\be \label{errorU}
e^{-i t H /\eps}-U^\eps_t=O\left( \frac{t^3}{\eps^3}\right),
\ee
so that $\tilde \varsigma^{\varepsilon}_{1,n}$ is of order 
$e^{-\delta t/\eps^2} (\delta t)^3/\eps^3$. Hence, $\eps^{-1 }\tilde \varsigma^{\varepsilon}_{1,n}$ is of order $e^{-\delta t/2\eps^2} \delta t$, which is an $O(\eps^2 \wedge \delta t )$ just as the current of $\varsigma^{\varepsilon}(t_n,t_{n+1})$ in the previous Section. The term $R_2$ therefore introduces an overall error on $\mn^{\varepsilon}_{n}-\tilde \mn^\eps$ of the same order as that of $\eps^{-1} \nabla \cdot j[S_{\varepsilon,(r-t_n)/\varepsilon^2}[\boldsymbol{\theta}_{n}^{\varepsilon}]]$ in the previous Section, that is an $O( \eps^2 \wedge \delta t)$. 


According to \fref{errorU}, the term $R_3$ is bounded by, using that $a_\eps=O(\eps^2 \wedge \delta t)$,
$$
\frac{a_\eps^3}{\eps^4}\int_{t_n}^{t_{n+1}}e^{-\frac{r-t_n}{\varepsilon^2}} dr \lesssim \eps^2 \wedge \delta t^2,
$$
which after summation from $n=0$ to $n=N-1$ gives a term of order $\eps^2 \wedge \delta t$. The error introduced by $R_3$ on $\mn^{\varepsilon}_{n}-\tilde \mn^\eps$ is therefore of order  $O( \eps^2 \wedge \delta t)$.



Regarding $R_4$, 
assuming $\mathsf{A}[\tilde \mn^{\varepsilon}(t)]$ and $\mn^{\varepsilon}(t)$ are differentiable w.r.t. $t$ with uniform bounds, the local error is an $O(\delta t^2)$, resulting in a global error on $\mn^{\varepsilon}_{n}-\tilde \mn^\eps$  of order $O(\delta t)$. 

Collecting previous estimates, we find that, for $0 \leq n \leq N$,
$$
\|\tilde \mn^\eps(t_n) - \mn^\eps_n \|_{L^1} \lesssim \delta t,
$$
which combined with \fref{interm}, gives
$$
\|\mn^\eps(t_n) - \mn^\eps_n \|_{L^1} \lesssim \delta t.
$$
In terms of the density operators, this translates to
$$
\|\tilde \varrho^{\varepsilon}(t_{n})-\varrho^{\varepsilon}_{n} \|_{\calE^2} \lesssim \delta t,
$$
which together with \fref{intermthet}, yields
$$
\|\varrho^{\varepsilon}(t_{n})-\varrho^{\varepsilon}_{n} \|_{\calE^2} \lesssim \delta t.
$$
This concludes our proof of Theorem \ref{th}.

\section{Conclusion} We have derived an asymptotic preserving scheme for the quantum Liouville-BGK equation in the diffusion limit. Our approach is quite general and applies to domains of arbitrary dimensions, and the key ingredient is the derivation of an equation for the local density. For this, we follow the asymptotic expansions obtained in the diffusion limit, which eventually lead to a modified quantum-drift diffusion equation. Once the density is known, it is injected into an asymptotic expansion of the density operator that is uniformly accurate in the diffusion parameter. The cost of the method for one time iteration is comparable to that of a split-step algorithm.

We have moreover derived theoretical error estimates under reasonable assumptions on the well-posedness of the different involved PDEs. The numerical implementation of the AP scheme will be addressed in a forthcoming work.

\bibliographystyle{plain}
\bibliography{bibliography}
 \end{document}